\documentclass[10.5 pt, a4paper]{amsart}
\usepackage{amssymb,latexsym,amsmath,amsfonts,amsthm, comment}

\usepackage{tikz}
\usepackage{multicol}
\usepackage{todonotes}
\usetikzlibrary{arrows}
\usetikzlibrary{decorations.markings}

\newcommand{\cev}[1]{\reflectbox{\ensuremath{\vec{\reflectbox{\ensuremath{#1}}}}}}

\numberwithin{equation}{section}

\usepackage[mathscr]{eucal}

\usepackage[T1]{fontenc}
\usepackage{txfonts}
\usepackage[pdftex,bookmarks=true]{hyperref}
\usepackage{graphicx}
 \usepackage{float}
\usepackage{enumerate}

\usepackage[all]{xy}

\usepackage{enumerate}

\theoremstyle{plain}
\newtheorem{theorem}{Theorem}[section]

\newtheorem{lemma}[theorem]{Lemma}

\theoremstyle{remark}

\newtheorem{exa}[theorem]{Example}

\theoremstyle{definition}
\newtheorem{dfn}[theorem]{Definition}

\author{RAKESH KUMAR}
\address{
School of Mathematics \& Computer Science\\
Indian Institute of Technology Goa\\
At Goa College of Engineering Campus\\
Farmagudi, Ponda-403401 \\
Goa, India} 
\email{rakesh20232101@iitgoa.ac.in}
\author{SHIV PARSAD}
\address{
School of Mathematics \& Computer Science\\
Indian Institute of Technology Goa\\
At Goa College of Engineering Campus\\
Farmagudi, Ponda-403401 \\
Goa, India} 

\email{shiv@iitgoa.ac.in}

\begin{document}
\title{filling systems of maximum size }

\subjclass[2000]{Primary 57M15; Secondary 05C10}

\keywords{Surface, filling system, fat graph}

\begin{abstract}
Let $S_g$ be a closed orientable surface of genus $g\geq 2$. A collection $\Omega = \{ \gamma_1, \dots, \gamma_s\}$ of pairwise non-homotopic simple closed curves on $S_g$ such that $\gamma_i$ and $\gamma_j$  are in minimal position, is called a \emph{filling system} or a \emph{filling} of $S_g$ if the complement $S_g\setminus \Omega$ is a disjoint union of $b$ topological discs for some $b\geq 1$. The \emph{size} of a filling system is defined as the number of its elements. We prove that the maximum size of a filling system on $S_g$ with $ 1 \leq b \leq 2g-2$ boundary components is $2g+b-1$. 
Furthermore, we give a lower bound on mapping class group orbits of filling systems of maximum size with $ 1 \leq b \leq g-2$  boundary components. 
\end{abstract}

\maketitle


\section{Introduction}
Let $S_g$ be a closed orientable surface of genus $g$. A set $\Omega = \{ \gamma_1, \dots, \gamma_s\}$ of pairwise non-homotopic simple closed curves on $S_g$ is called a \emph{filling system} or a \emph{filling} of $S_g$, if the complement $S_g\setminus \Omega$ is a disjoint union of $b$ topological discs for some $b\geq 1$. It is assumed that the curves in a filling system are pairwise in minimal position, i.e. $i(\gamma_i, \gamma_j )= \lvert \gamma_i \cap \gamma_j \rvert$ for $i \neq j$. We call it a minimal filling if the complement is a single disc. Furthermore, if $\Omega$ contains exactly two curves, it is called a filling pair.
Filling systems on closed surfaces have increasingly gained significance in the study of the mapping class group of surfaces and the moduli space of hyperbolic surfaces via the systolic function. The study of filling systems originates from the work of Thurston~\cite{WT}, where the author has defined a subset  of Teichmüller space
 space using filling systems which he claimed to be a spine of the Teichm\"ller space~\cite{II}. In~\cite{TW2}, Thurston used filling pairs to construct pseudo-Anosov mapping classes~\cite{FM}. This construction was generalized by Penner~\cite{RCP} to filling systems of any size. Recently, filling systems on surfaces has been studied by Schaller~\cite{Schmutz}, Aougab-Huang~\cite{TA}, Fanoni-Parlier~\cite{Parlier}, Sanki~\cite{BS}, Parsad-Sanki \cite{PS}, Bourque \cite{BM} and others. 

 In an unpublished version \cite{parsad2017filling}, the authors provide an upper bound on the maximum size of a filling system with given boundary components and also proved the sharpness of the upper bound. But there is a gap in proof of the sharpness of upper bound. In the following result, we prove that under some restriction on number of boundary components, the upper bound is sharp. 

\begin{theorem}\label{Main}
Let $S_g$ be a closed orientable surface of genus $g\geq2$ and let $U_{g,b}$ denotes the maximum size of filling systems on $S_{g}$ with $b$ boundary components. Then $U_{g,b}=2g+b-1$, for $ 1 \leq b \leq 2g-2$.
\end{theorem}

 There is an action of mapping class group on the set of filling systems on $S_g$ of fixed size and given number of boundary components. In \cite{TA}, the authors prove that there are at least exponentially many mapping class group orbits of minimal filling pairs. In the following result, we provide a lower bound on mapping class group orbits of filling systems of maximum size with given boundary components.  
\begin{theorem} \label{orbit}
Let $S_g$ be a closed orientable surface of genus $g\geq2$ and $N_{g,b}$  denotes the mapping class orbits of filling systems of $S_g$ of size $2g+b-1$ with $b$ boundary components. Then:\begin{enumerate}
 \item $N_{g,1} \geq |\Re^{s}_{g-s}| + 1$, where $\Re^{s}_{g-s}$ is the union of equivalence classes of integer partitions of $ g - s $ into $s$ non-negative integers,  with equivalence under permutation by $ S_s $, where $2\leq s \leq g$.
\item $N_{g,2} \geq |\Re^{s}_{g-s+1}| + 1 $, where $\Re^{s}_{g-s+1}$ is the union of equivalence classes of integer partitions of $ g - s +1$ into $s$ non-negative integers, with equivalence under permutation by $ S_s $, where $2\leq s \leq g+1$. 
\item $N_{g,b} \geq |\Re^{b}_{g-1}| + 1$, where $ \Re^{b}_{g-1} $ is the union of equivalence classes of integer partitions of $g-1$ into $b$ non-negative integers for $3 \leq b \leq g-2$, with equivalence under permutation by $ S_b $.
\end{enumerate}

\end{theorem}

\section{Fat graphs}
In this section, we recall some definitions from topological graph theory, in particular fat graphs. 

\begin{dfn}
A fat  graph is a quadruple $\Gamma=(E, \sim, \sigma_1, \sigma_{0})$, where 
\begin{enumerate}
\item $E=\{\vec{e}_1, \cev{e}_1, \dots, \vec{e}_n, \cev{e}_n\}$ is a set of even cardinality called the set of directed edges,

\item  $\sim$ is an equivalence relation on $E$, 

\item  $\sigma_1:E\to E$ is a fixed point free involution sending $\vec{e}_i $ to $\cev{e}_i$ for each $i$. We think of fixed point free involution as reversing the direction of the directed edge.
\item $\sigma_{0}$  is a permutation on $E$, whose cycles correspond to cyclic order on the set of directed edges emanating from each vertex. 
\end{enumerate}
The set $V:=E/\!\!\sim$ of equivalence classes of $\sim$ is the set of vertices. We note that each equivalence class $v\in V$ is the set consisting of all directed edges emanating from the corresponding vertex $v$. The set $E_1:=E/\sigma_1$ is the set of undirected edges of $\Gamma$. The \emph{degree} of a vertex $v$ is defined as its cardinality.
  
\end{dfn}
A fat graph is called \emph{decorated} if the degree of each vertex is an even integer and at least $4$. A cycle in a decorated fat graph is called \emph{standard} if every two consecutive edges in the cycle are opposite to each other in the cyclic order at their shared vertex. Given a fat graph $\Gamma$, the set of boundary components is denoted by $\partial\Gamma$, which corresponds to the set of cycles of the permutation $ \sigma_{\infty}=\sigma_0*\sigma_1$ (see Lemma 2.4 in~\cite{BS}).\\
We can construct a unique topological surface with boundary corresponding to a fat graph by thickening the edges of a fat graph $\Gamma$. Thus, we can talk about the number of boundary components, genus, and other topological properties. For more details on fat graphs, we refer the reader to ~\cite[Section 2]{BS}, \cite{JS}, \cite{JS1} and \cite{MP}.

\begin{exa}\label{eg:triple_2}
Consider the  fat graph $\Gamma=(E,\sim,\sigma_1,\sigma_0)$ of genus $2$ (see Figure~\ref{eg:1}), described as follows:

\begin{enumerate}
\item $E=\{\vec{e}_i, \cev{e}_i, \vec{f}_i, \cev{f}_i\;|\;i=1,2,3\}$.
\item The set of equivalence classes of $\sim$ is $V=\{v_1, v_2, v_3\}$, where 
$ v_1=\left\{\vec{e}_1, \cev{f}_2, \cev{e}_3, \vec{f}_1\right\},\\\,v_2=\left\{\vec{e}_2, \cev{f}_1, \cev{e}_1, \vec{f}_2\right\},\text{ and }v_3=\left\{\vec{e}_3, \cev{f}_3, \cev{e}_2, \vec{f}_3\right\}.$
\item The fixed point free involution $\sigma_1$ is defined by $\sigma_1(\vec{e}_i)=\cev{e}_i$ and $\sigma_1(\vec{f}_i)=\cev{f}_i$ for $i=1,2,3.$
\item The permutation (fat graph structure) $\sigma_0$ is given by $ \sigma_0=\prod\limits_{i=1}^{3}\sigma_{v_i}$, where 
$ \sigma_{v_1}=\left(\vec{e}_1, \cev{f}_2, \cev{e}_3, \vec{f}_1\right),\\\,\sigma_{v_2}=\left(\vec{e}_2, \cev{f}_1, \cev{e}_1, \vec{f}_2\right),\text{ and }\sigma_{v_3}=\left(\vec{e}_3, \cev{f}_3, \cev{e}_2, \vec{f}_3\right).$
\end{enumerate}
It is straightforward to see from Figure~\ref{eg:1} that $\Gamma_0$ has one boundary component and $3$ standard cycles of lengths $3,2$ and $1$. Therefore, $\Gamma_0$ corresponds to a minimal filling triple of $S_2$.
\begin{figure}[htbp]
    \centering
    \includegraphics[width=0.5\linewidth]{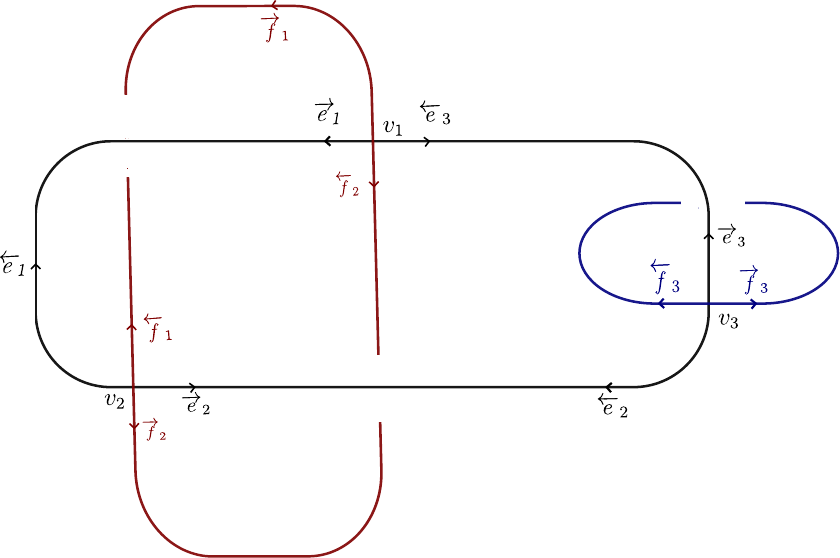}
    \caption{$\Gamma_0$}
    \label{eg:1}
\end{figure}
\end{exa}
\section{Filling Systems Of Maximum Size}
In this section, we prove Theorem \ref{Main}. We need the following lemma:
\begin{lemma}
If $U_{g,b}$ denotes the maximum size of filling systems on $S_{g}$ with $b$ boundary components, then $U_{g,b} \leq 2g+b-1$, for all $g,b \in \mathbb{N}$.
 \end{lemma}
\noindent For proof see \cite[ Lemma 3.1]{parsad2017filling}.

In order to prove Theorem \ref{Main} we need to construct a 4-regular fat graph of genus $g$ with $2\leq b\leq 2g-2$ boundary components and $2g+b-1$ standard cycles. 
\begin{proof} The sharpness of bound $U_{g,1}$ realized by a chain of order $2g$ and the sharpness of bound $U_{g,2}$ realized by a chain of order $2g+1$. Now, for $ 3 \leq b \leq 2g-2$, we will prove this in two cases.\\ 
{\bf Case 1}:  $g-1 \leq b \leq 2g-2.$\\
First, we consider the graph for $\Gamma_{(g,g-1)} = (E, \sim, \sigma_1, \sigma_0)$ for $b=g-1$ described below (see Figure 2 ). 
 \begin{enumerate}
\item  $E=\{ \vec{e}_i, \cev{e}_i  \vec{f}_j, \cev{f}_j \rvert 1\leq i \leq 2g-2 \text{ and } 1 \leq j \leq 4g-4 \}$ is the set of directed edges.
 \item $\sim$ is determined by the partition $V= \{ v_i | 1 \leq i \leq 3g-3 \}$, where 
 \begin{align*}
 v_1 &= \{ \vec{e}_1 , \cev{f}_1, \cev{e}_{2g-2}, \vec{f}_1 \}, \\
v_{2i} & = \{ \vec{e}_{2i} , \cev{f}_{4i-1}, \cev{e}_{2i-1}, \vec{f}_{4i-2} \}; 1\leq i  \leq g-1,  \\
v_{2i+1} & = \{ \vec{e}_{2i+1} , \cev{f}_{4i+1}, \cev{e}_{2i}, \vec{f}_{4i+1} \}; 1\leq i  \leq g-2 
 \text{ and }\\
v_{2g-2+i} & = \{ \cev{f}_{4i-2} , \vec{f}_{4i}, \vec{f}_{4i-1}, \cev{f}_{4i} \}; 1\leq i  \leq g-1  
\end{align*}
\item  $\sigma_1(\vec{e}_{i})= \cev{e}_{i} $ and $\sigma_1(\vec{f}_{j})= \cev{f}_{j}; 1\leq i \leq 2g-2 \text{ and } 1 \leq j \leq 4g-4,  $
\item $\sigma_{0}=C_1 \dots C_{3g-3,}$ where
\begin{align*}
C_1 &= (\vec{e}_1 , \cev{f}_1, \cev{e}_{2g-2}, \vec{f}_1  ), \\ 
C_{2i} & = (\vec{e}_{2i} , \cev{f}_{4i-1}, \cev{e}_{2i-1}, \vec{f}_{4i-2} ); 1\leq i  \leq g-1,\\  
C_{2i+1} &= ( \vec{e}_{2i+1} , \cev{f}_{4i+1}, \cev{e}_{2i}, \vec{f}_{4i+1} ); 1\leq i  \leq g-2, and\\ 
C_{2g-2+i} & = (\cev{f}_{4i-2} , \vec{f}_{4i}, \vec{f}_{4i-1}, \cev{f}_{4i} ); 1\leq i  \leq g-1.
\end{align*}
\end{enumerate}
 \begin{figure}[htbp]
     \centering
     \includegraphics[width=0.6\linewidth]{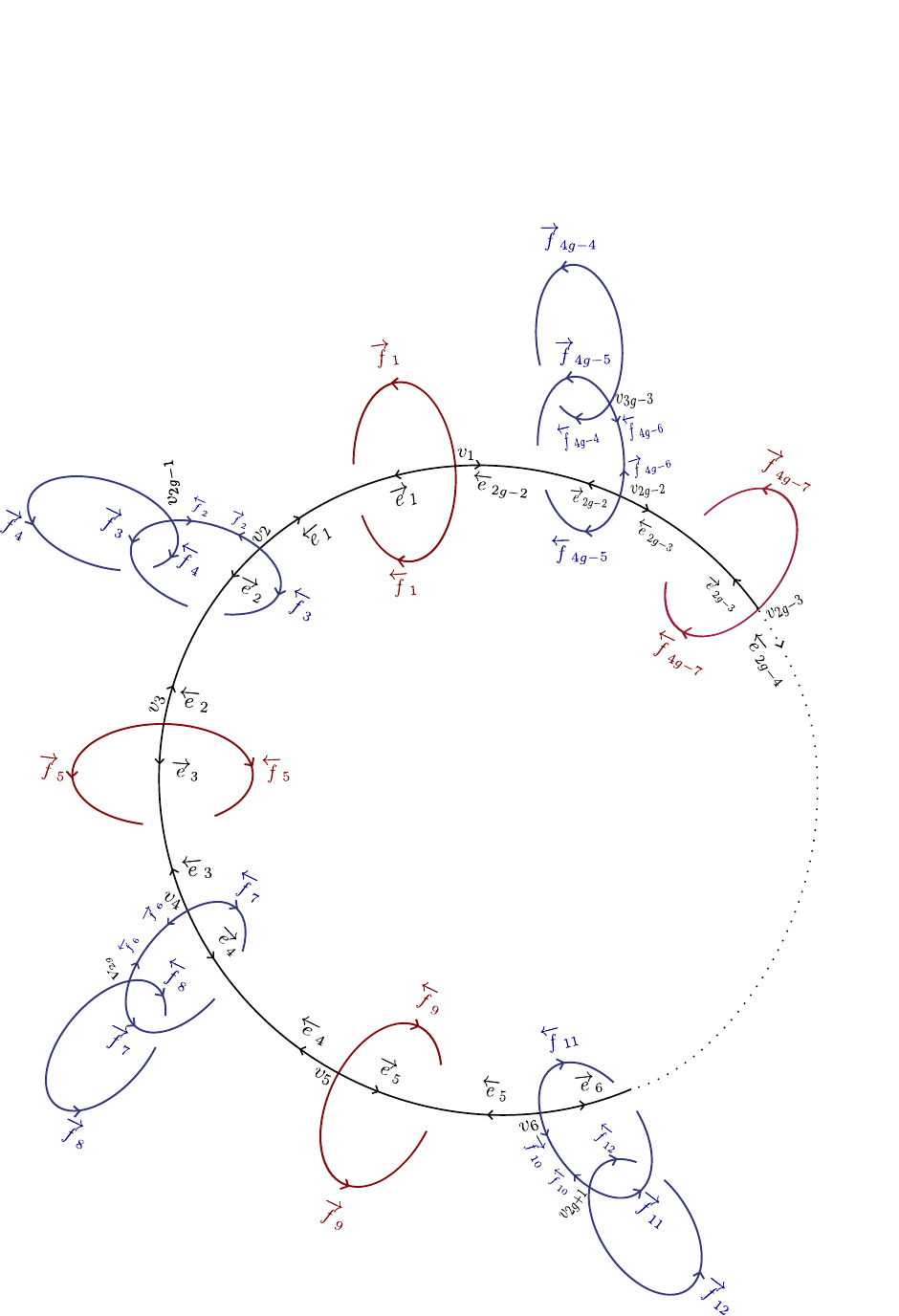}
     \caption{The fat graph $\Gamma_{(g,g-1)}$}
     
 \end{figure}
Now, we will count the boundary components of $\Gamma_{(g,g-1)}$ using the permutation $\sigma_{\infty}$. The boundary components of the graph $\Gamma_{(g,g-1)}$ given by: \\
$\partial_i= \vec{e}_{2i-1} \vec{f}_{4i-2} \vec{f}_{4i} \cev{f}_{4i-2} \vec{e}_{2i} \vec{f}_{4i+1} \cev{e}_{2i} \cev{f}_{4i-1} \cev{f}_{4i} \vec{f}_{4i-1} \cev{e}_{2i-1} \cev{f}_{4i-3} $ where $1 \leq i \leq g-1$.\\
So, there are $g-1$ boundary components for graph $\Gamma_{1}$.

Now we will give a fat graph for $b=g-1+k$ where $1 \leq k \leq g-1$.
we consider the graph for $\Gamma_{(g,b)} = (E, \sim, \sigma_1, \sigma_0)$ for $g\leq b \leq 2g-2$ described below (see Figure 3 ). 
\begin{enumerate}
\item  $E=\{ \vec{e}_i, \cev{e}_i  \vec{f}_j, \cev{f}_j \rvert 1\leq i \leq 2g-2 \text{ and } 1 \leq j \leq 4g+2k-4 \}$ is the set of directed edges.
\item $\sim$ is determined by the partition $V= \{ v_i | 1 \leq i \leq 3g-3+k \}$, where
\begin{align*}
v_1 & = \{ \vec{e}_1 , \cev{f}_1, \cev{e}_{2g-2}, \vec{f}_1 \}, \\
v_{2k-1} & = \{ \vec{e}_{2k-1} , \cev{f}_{6k-5}, \cev{e}_{2k-2}, \vec{f}_{6k-5} \}; 2 \leq k \leq g-1, \\
v_{2k} & = \{ \vec{e}_{2k} , \cev{f}_{6k+1}, \cev{e}_{2k-1}, \vec{f}_{6k+1} \}; 1 \leq k \leq g-1   ,\\
v_{2k+2l-1} & = \{ \vec{e}_{2k+2l-1} , \cev{f}_{6k+4l-3}, \cev{e}_{2k+2l-2}, \vec{f}_{6k+4l-3} \}; 1\leq l  \leq g-1-k ,\\ 
v_{2k+2l} & = \{ \vec{e}_{2k+2l} , \cev{f}_{6k+4l-1}, \cev{e}_{2k+2l-1}, \vec{f}_{6k+4l-2} \}; 1\leq l  \leq g-1-k ,\\ 
v_{2g-1} & = \{ \cev{f}_{2} , \vec{f}_{4}, \vec{f}_{3}, \cev{f}_{5} \}\\
v_{2g-2+2m} & = \{ \cev{f}_{6m-2} , \vec{f}_{6m}, \vec{f}_{6m-1}, \cev{f}_{6m}\};1\leq m  \leq k , \\
v_{2g+2m-1} & =\{ \cev{f}_{6m+2} , \vec{f}_{6m+4}, \vec{f}_{6m+3}, \cev{f}_{6m+5}\};1\leq m  \leq k-1, \text{ and } \\
v_{2g-2+2k+n}&   = \{ \cev{f}_{6k+4n-2} , \vec{f}_{6k+4n}, \vec{f}_{6k+4n-1}, \cev{f}_{6k+4m}\};1\leq n  \leq g-1-k.
\end{align*}
\item  $\sigma_1(\vec{e}_{i})= \cev{e}_{i} $ and $\sigma_1(\vec{f}_{j})= \cev{f}_{j}; 1\leq i \leq 2g-2 \text{ and } 1 \leq j \leq 4g+2k-4.$
\item $\sigma_{0}=C_1 \dots C_{3g-3+k,}$ where
\end{enumerate}
\begin{align*}
C_1 & = (\vec{e}_1 , \cev{f}_1, \cev{e}_{2g-2}, \vec{f}_1 ), \\
C_{2k-1} &= ( \vec{e}_{2k-1} , \cev{f}_{6k-5}, \cev{e}_{2k-2}, \vec{f}_{6k-5} ) ; 2 \leq k \leq g-1, \\
C_{2k} & = (  \vec{e}_{2k} , \cev{f}_{6k+1}, \cev{e}_{2k-1}, \vec{f}_{6k+1}  ); 1 \leq k \leq g-1,\\  
C_{2k+2l-1}& =(\vec{e}_{2k+2l-1} , \cev{f}_{6k+4l-3}, \cev{e}_{2k+2l-2}, \vec{f}_{6k+4l-3} ); 1\leq l  \leq g-1-k , \\ 
C_{2k+2l}&=(\vec{e}_{2k+2l} , \cev{f}_{6k+4l-1}, \cev{e}_{2k+2l-1}, \vec{f}_{6k+4l-2} ); 1\leq l  \leq g-1-k , \\
C_{2g-1}&= (\cev{f}_{2} , \vec{f}_{4}, \vec{f}_{3}, \cev{f}_{5} ),\\
C_{2g-2+2m}&= (\cev{f}_{6m-2} , \vec{f}_{6m}, \vec{f}_{6m-1}, \cev{f}_{6m});1\leq m  \leq k ,\\
C_{2g+2m-1}&=( \cev{f}_{6m+2} , \vec{f}_{6m+4}, \vec{f}_{6m+3}, \cev{f}_{6m+5} );1\leq m  \leq k-1, \text{ and } \\
C_{2g-2+2k+n}&= (\cev{f}_{6k+4n-2} , \vec{f}_{6k+4n}, \vec{f}_{6k+4n-1}, \cev{f}_{6k+4m});1\leq n  \leq g-1-k.
\end{align*}
\begin{figure}[htbp]
    \centering
    \includegraphics[width=0.98\linewidth]{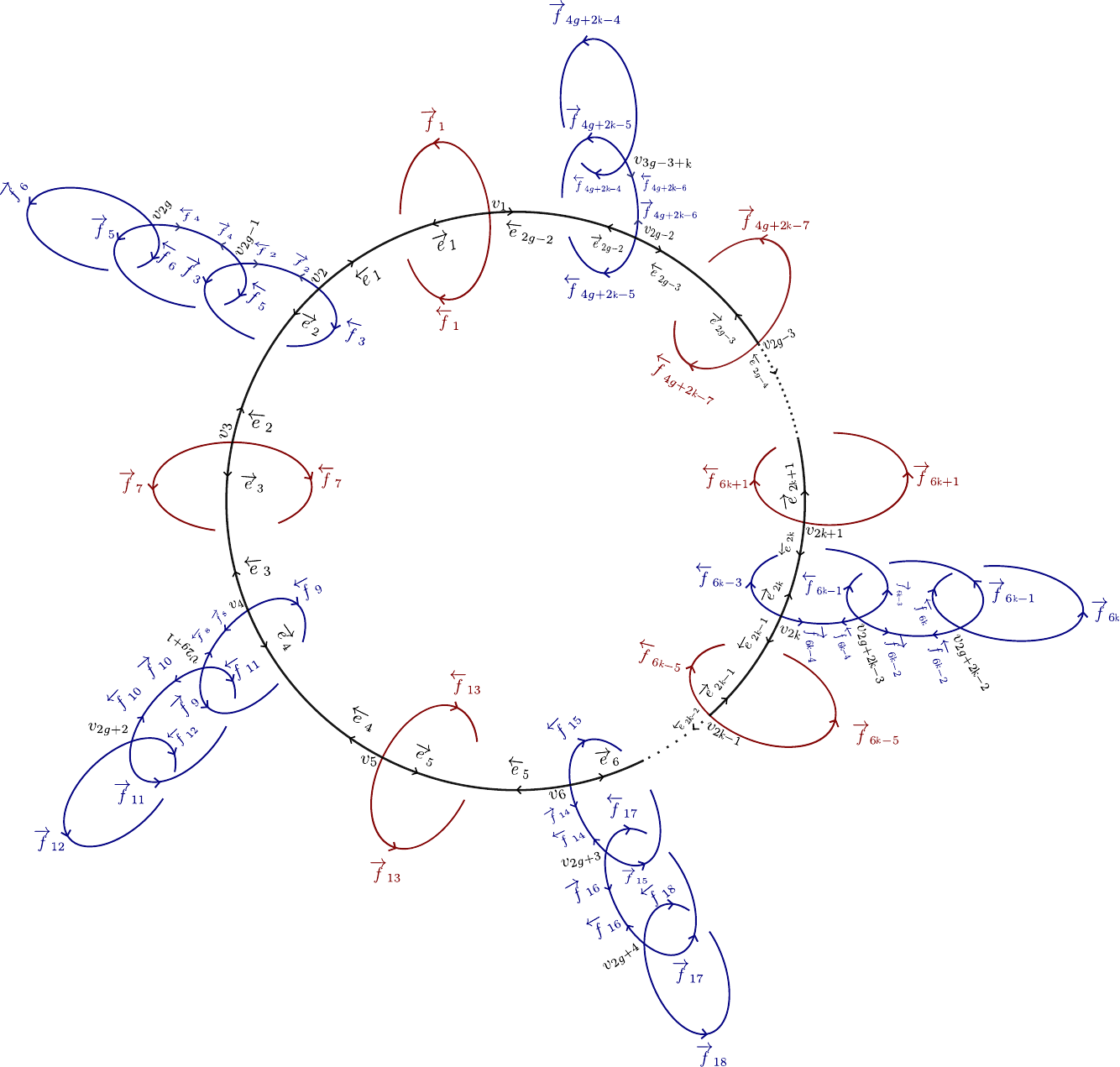}
    \caption{The fat graph $\Gamma_{(g,b)}$ \text{ for }$g\leq b \leq 2g-2$}
    
\end{figure}

Now, we count the boundary components of $\Gamma_{(g,b)}$ using the permutation $\sigma_{\infty}$. We will consider two subcases.\\
{\bf Subcase 1}:  $1 \leq k \leq g-2.$\\
The boundary components of the graph $\Gamma_{(g,b)}$ given by: 
\begin{align*}
\partial_{1}& = \vec{e}_1 \vec{f}_2 \vec{f}_4 \vec{f}_6 \cev{f}_4 \vec{f}_3 \cev{e}_1 \cev{f}_1, \\
\partial_2 &= \vec{e}_2 \vec{f}_7 \cev{e}_2 \cev{f} _3 \cev{f}_5 \cev{f}_6 \vec{f}_5 \cev{f}_2, \\
 \vdots  \\
\partial_{2k-1} &= \vec{e}_{2k-1} \vec{f}_{6k-4} \vec{f}_{6k-2} \vec{f}_{6k} \cev{f}_{6k-2} \vec{f}_{6k-3} \cev{e}_{2k-1} \cev{f}_{6k-5},\\
\partial_{2k} & = \vec{e}_{2k} \vec{f}_{6k+1} \cev{e}_{2k} \cev{f}_{6k-3} \cev{f}_{6k-1} \cev{f}_{6k} \vec{f}_{6k-1} \cev{f}_{6k-4},\\
\partial_{2k+l} & = \vec{e}_{2k+2l-1} \vec{f}_{6k+4l-2} \vec{f}_{6k+4l} \cev{f}_{6k+4l-2} \vec{e}_{2k+2l} \vec{f}_{6k+4l+1} \cev{e}_{2k+2l} \cev{f}_{6k+4l-1} \cev{f}_{6k+4l} \vec{f}_{6k+4l-1} \cev{e}_{2k+2l-1} \cev{f}_{6k+4j-3} \text{ where } 1 \leq l \leq g-1-k.
\end{align*}
{\bf Subcase 2}:  $ k = g-1.$\\
The boundary components of the graph $\Gamma_{(g,b)}$ given by: 
\begin{align*}
\partial_{1}& = \vec{e}_1 \vec{f}_2 \vec{f}_4 \vec{f}_6 \cev{f}_4 \vec{f}_3 \cev{e}_1 \cev{f}_1, \\
\partial_2 &= \vec{e}_2 \vec{f}_7 \cev{e}_2 \cev{f} _3 \cev{f}_5 \cev{f}_6 \vec{f}_5 \cev{f}_2, \\
 \vdots  \\
\partial_{2k-1} &= \vec{e}_{2k-1} \vec{f}_{6k-4} \vec{f}_{6k-2} \vec{f}_{6k} \cev{f}_{6k-2} \vec{f}_{6k-3} \cev{e}_{2k-1} \cev{f}_{6k-5},\\
\partial_{2k} & = \vec{e}_{2k} \vec{f}_{6k+1} \cev{e}_{2k} \cev{f}_{6k-3} \cev{f}_{6k-1} \cev{f}_{6k} \vec{f}_{6k-1} \cev{f}_{6k-4}.
\end{align*}

So, there are $b=2k+g-1-k=g-1+k$ boundary components for graph $\Gamma_{(g,b)}$.

{\bf Case 2.} $3 \leq b \leq g-2.$

Consider the graph $\Gamma_{(g,b)} = (E, \sim, \sigma_1, \sigma_0)$ described below (see Figure 4 ). 
\begin{figure}[htbp]
     \centering
     \includegraphics[width=1.2\linewidth]{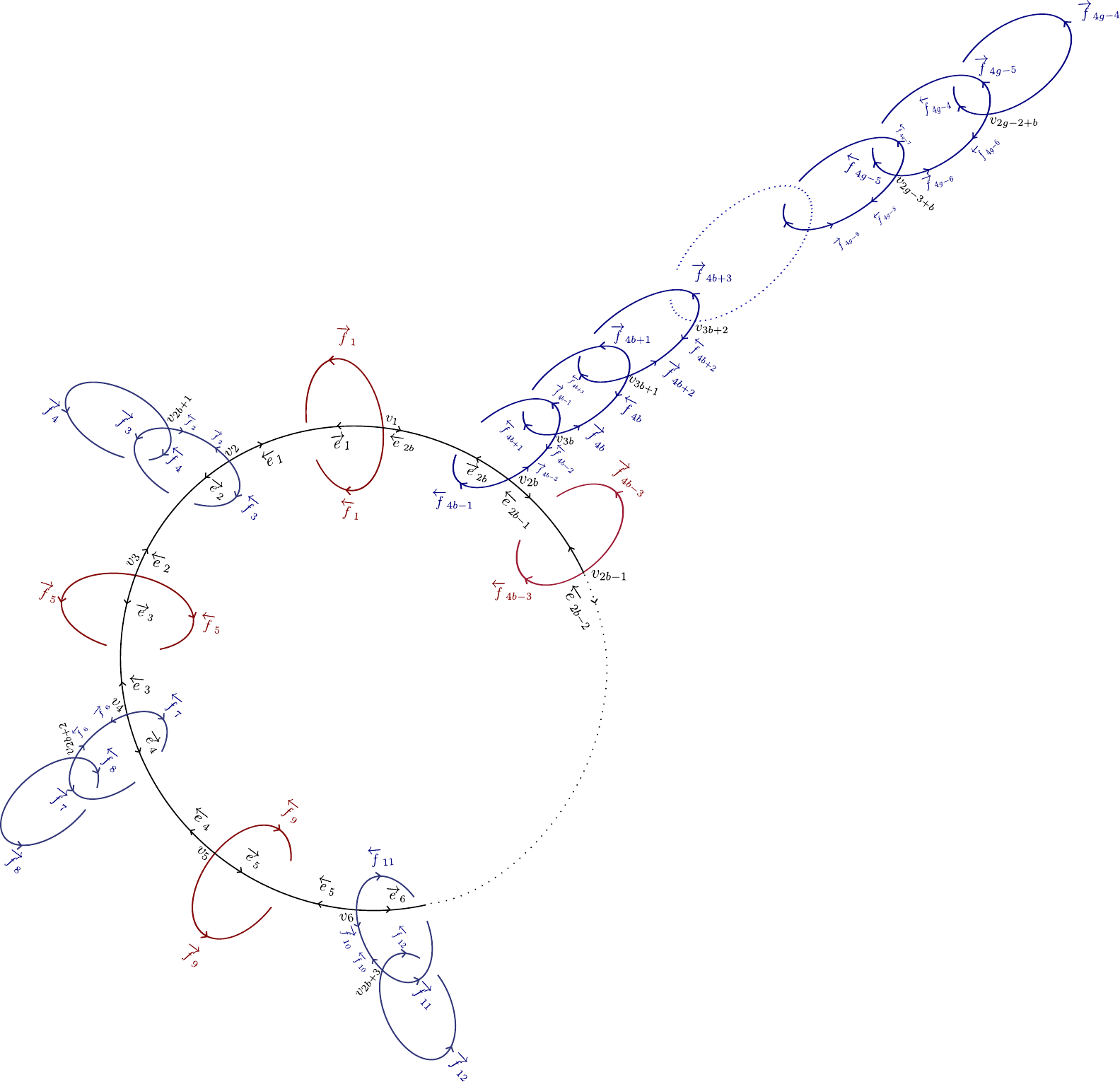}
     \caption{The fat graph $\Gamma_{(g,b)}$ \text{ for } $3 \leq b \leq g-2$}
     
 \end{figure}
 \begin{enumerate}
 \item  $E=\{ \vec{e}_i, \cev{e}_i  \vec{f}_j, \cev{f}_j \rvert 1 \leq i \leq 2g-4 \text{ and } 1 \leq j \leq 4g-4 \}$ is the set of directed edges.
\item $\sim$ is determined by the partition $V= \{ v_i | 1 \leq i \leq 2g-2+b \}$, where
\begin{align*}
v_1 &= \{ \vec{e}_1 , \cev{f}_1, \cev{e}_{2b}, \vec{f}_1 \}, \\
v_{2i-1} &= \{ \vec{e}_{2i-1} , \cev{f}_{4i-3}, \cev{e}_{2i-2}, \vec{f}_{4i-3} \}; 1\leq i  \leq b,  \\
v_{2i} &= \{ \vec{e}_{2i} , \cev{f}_{4i-1}, \cev{e}_{2i-1}, \vec{f}_{4i-2} \}; 1\leq i  \leq b,  \\
v_{2b+j} &= \{ \cev{f}_{4j-2} , \vec{f}_{4j}, \vec{f}_{4j-1}, \vec{f}_{4j} \}; 1\leq j  \leq b-1, \\
v_{3b-1+l} &= \{ \cev{f}_{4b+2l-4} , \vec{f}_{4i+2b-2}, \vec{f}_{4b+2l-3}, \cev{f}_{4b+2l-1} \}; 1\leq l  \leq 2g-2b-2, \text{ and }  \\
v_{2g+b-2} &= \{ \cev{f}_{4g-6} , \vec{f}_{4g-4}, \vec{f}_{4g-5}, \cev{f}_{4g-4} \}. 
 \end{align*}
\item  $\sigma_1(\vec{e}_{i})= \cev{e}_{i} $ and $\sigma_1(\vec{f}_{j})= \cev{f}_{j}; 1\leq i \leq 2g-4 \text{ and } 1 \leq j \leq 4g-4 .$
\item $\sigma_{0}=C_1 \dots C_{2g+i-2,}$ where 
\end{enumerate}
\begin{align*}
C_1 &= ( \vec{e}_1 , \cev{f}_1, \cev{e}_{2b}, \vec{f}_1 ), \\ 
C_{2i-1} &= ( \vec{e}_{2i-1} , \cev{f}_{4i-3}, \cev{e}_{2i-2}, \vec{f}_{4i-3}  ); 1\leq i  \leq b,\\  
C_{2i} &= (\vec{e}_{2i} , \cev{f}_{4i-1}, \cev{e}_{2i-1}, \vec{f}_{4i-2}); 1\leq i  \leq b,  \\
C_{2b+j} &= ( \cev{f}_{4j-2} , \vec{f}_{4j}, \vec{f}_{4j-1}, \vec{f}_{4j} ); 1\leq j  \leq b-1,  \\
C_{3b-1+l} &= ( \cev{f}_{4b+2l-4} , \vec{f}_{4i+2b-2}, \vec{f}_{4b+2l-3}, \cev{f}_{4b+2l-1}  ); 1\leq l  \leq 2g-2b-2,  \text{ and } \\
C_{2g+b-2} &= ( \cev{f}_{4g-6} , \vec{f}_{4g-4}, \vec{f}_{4g-5}, \cev{f}_{4g-4} ).
\end{align*}
Now, we count the boundary components of $\Gamma_{(g,b)}$ using the permutation $\sigma_{\infty}$. The boundary components of the graph $\Gamma_{(g,b)} \text{ for } 3 \leq b \leq g-2$ given by: \\
$\partial_i= \vec{e}_{2i-1} \vec{f}_{4i-2} \vec{f}_{4i} \cev{f}_{4i-2} \vec{e}_{2i} \vec{f}_{4i+1} \cev{e}_{2i} \cev{f}_{4i-1} \cev{f}_{4i} \vec{f}_{4i-1} \cev{e}_{2i-1} \cev{f}_{4i-3} $ where $1 \leq i \leq b-1 $.\\
$\partial_b= \vec{e}_{2b-1} \vec{f}_{4b-2} \vec{f}_{4b} \vec{f}_{4b+2} \vec{f}_{4b+4}  \dots \vec{f}_{4g-6} \vec{f}_{4g-4} \cev{f}_{4g-6} \vec{f}_{4g-7} \cev{f}_{4g-10} \vec{f}_{4g-11} \cev{f}_{4g-14} \vec{f}_{4g-15} \cev{f}_{4g-18} \dots \dots \\ \cev{f}_{4i-2} \vec{e}_{2b} \vec{f}_{1} \cev{e}_{2i} \cev{f}_{4i-1} \cev{f}_{4i+1} \cev{f}_{4i+3} \dots   \dots \cev{f}_{4g-5} \cev{f}_{4g-4} \vec{f}_{4g-5} \cev{f}_{4g-8} \vec{f}_{4g-9} \cev{f}_{4g-12} \vec{f}_{4g-13} \cev{f}_{4g-16} \vec{f}_{4g-17} \dots \vec{f}_{4b-1} \cev{e}_{2b-1} \cev{f}_{4b-3}.$\\
So, there are $3 \leq b \leq g-2$ boundary components for graph $\Gamma_{(g,b)}$.
\end{proof}
\noindent\textbf{Remark:} We observe that the condition $b\leq 2g-2$ can not be dropped from the statement of Theorem \ref{Main}. We will give an example of a filling system on $S_g$ with $b>2g-2$ boundary components and $U_{g,b} \neq 2g+b-1$. Let $g=2 \text{ and } b=3$. The number of vertices for $b=3$ and $g=2$ is 5 by the Euler characteristic argument. We will see all possible fat graphs with $6$ standard cycles, $5$ vertices and $10$ edges. Possible fat graphs for this combination are $\Gamma_1 ,\Gamma_2, \Gamma_3, \Gamma_4, \Gamma_{5}, \text{ and } \Gamma_{6}$(see figure 5, 6, 7, 8, 9 and 10 ). We can see that $\Gamma_1 ,\Gamma_2, \Gamma_3,  \text{ and } \Gamma_{4}$(see figure 5, 6, 7, and 8 ) are not filling for $S_2$ because curves are pairwise homotopic. Now, we can see $\Gamma_{6}$ is a chain of order $6$. As we know boundary components for an even order chain is $1$. Now we will count boundary components for $\Gamma_{5}$. 
Consider the graph for  $\Gamma_5 = (E, \sim, \sigma_1, \sigma_0)$ as described below (see figure 9 ). \\
\begin{figure}[htbp]
\begin{minipage}[t]{0.48\textwidth}
\includegraphics[width=\linewidth,keepaspectratio=true]{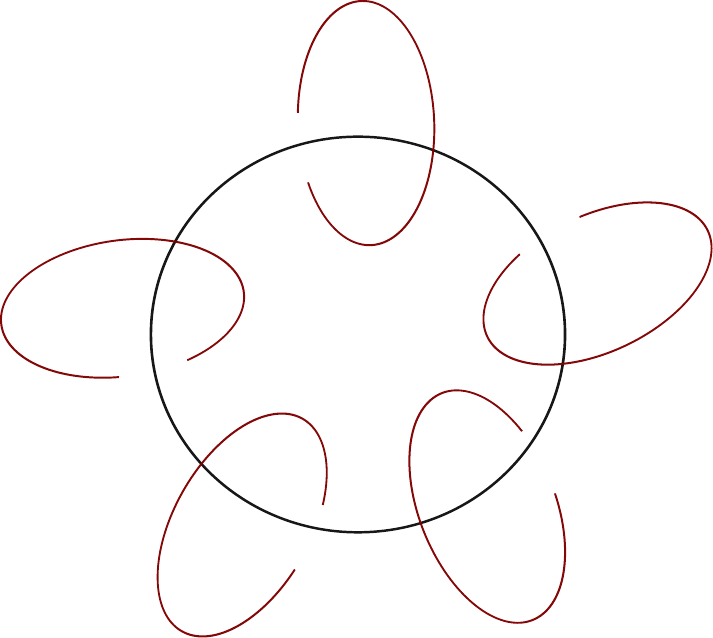}
\caption{$\Gamma_1$}

\end{minipage}
\hspace*{\fill} 
\begin{minipage}[t]{0.48\textwidth}
\includegraphics[width=\linewidth,keepaspectratio=true]{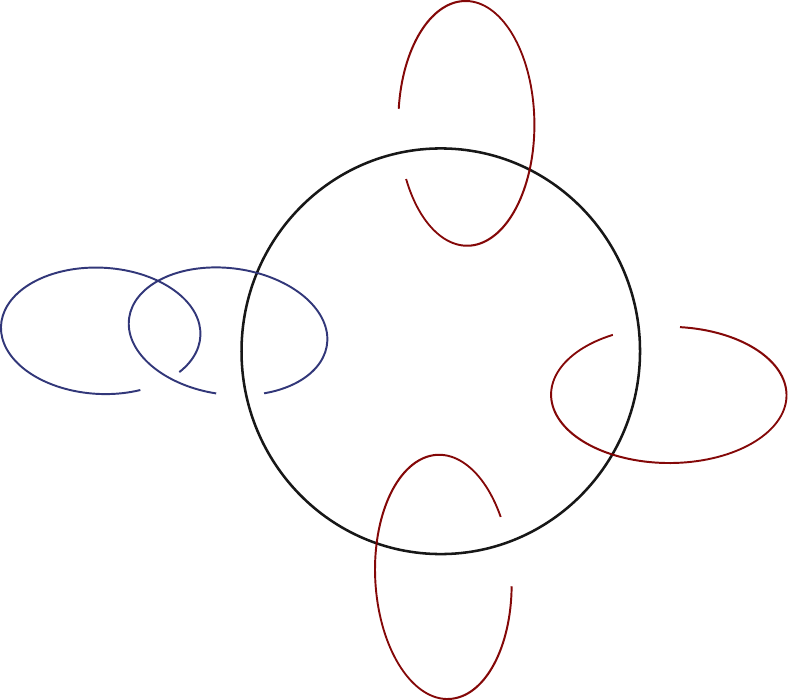}
\caption{$\Gamma_2$}

\end{minipage}
\end{figure}
\begin{figure}[htbp]
\begin{minipage}[t]{0.48\textwidth}
\includegraphics[width=\linewidth,keepaspectratio=true]{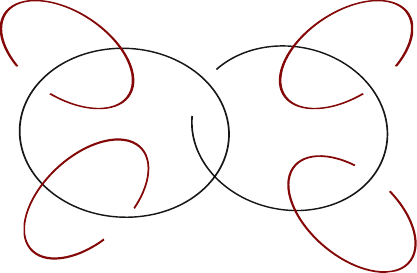}
\caption{$\Gamma_3$}

\end{minipage}
\hspace*{\fill} 
\begin{minipage}[t]{0.48\textwidth}
\includegraphics[width=\linewidth,keepaspectratio=true]{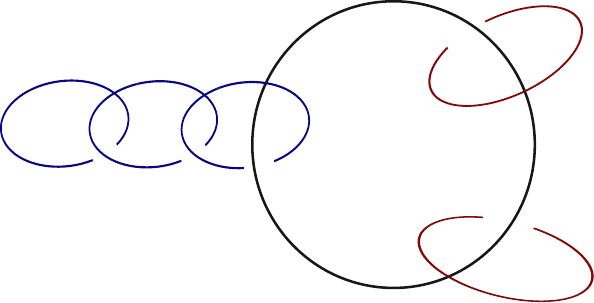}
\caption{$\Gamma_4$}

\end{minipage}
\end{figure}
\begin{figure}[htbp]
    \centering
    \includegraphics[width=0.7\linewidth]{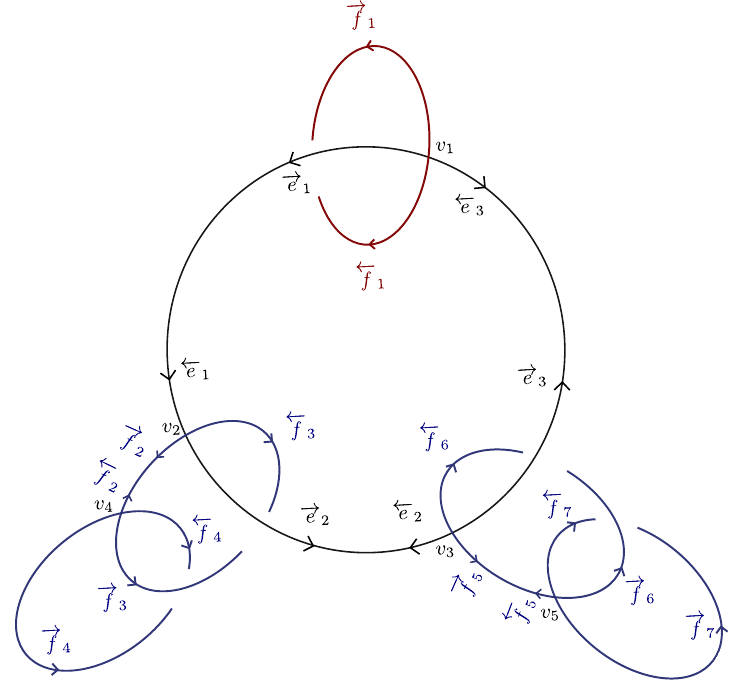}
    \caption{$\Gamma_5$}
    
\end{figure}
\begin{figure}[htbp]
    \centering
    \includegraphics[width=0.7\linewidth]{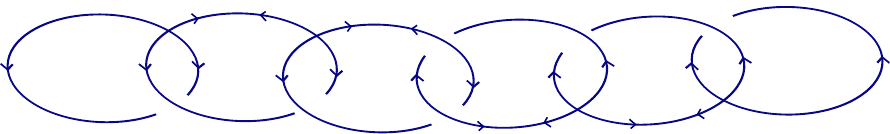}
    \caption{$\Gamma_6$}
    
\end{figure}

\begin{enumerate}
\item  $E=\{ \vec{e}_i, \cev{e}_i  \vec{f}_j, \cev{f}_j \rvert 1 \leq i \leq 3 \text{ and } 1 \leq j \leq 7  \}$ is the set of directed edges.
  
\item $\sim$ is determined by the partition $V= \{ v_i | 1 \leq i \leq 5 \}$, where
\begin{align*}
 v_1 &= \{ \vec{e}_1 , \cev{f}_1, \cev{e}_{3}, \vec{f}_1 \}, \\
v_{2} & = \{ \vec{e}_{2} , \cev{f}_{3}, \cev{e}_{1}, \vec{f}_{2} \}, \\
v_{3} & = \{ \vec{e}_{3} , \cev{f}_{6}, \cev{e}_{2}, \vec{f}_{5} \}, \\
v_{4} & = \{ \cev{f}_{2} , \vec{f}_{4}, \vec{f}_{3}, \cev{f}_{4} \}, \text{ and }\\
v_{5} & = \{ \cev{f}_{5} , \vec{f}_{7}, \vec{f}_{6}, \cev{f}_{7} \}. 
\end{align*}
\item  $\sigma_1(\vec{e}_{i})= \cev{e}_{i} $ and $\sigma_1(\overrightarrow f_{j})= \overleftarrow f_{j}; 1\leq i \leq 3 \text{ and } 1 \leq j \leq 7.$
\item $\sigma_{0}=C_1 \dots C_{5,}$ where 
\end{enumerate}
\begin{align*}
C_1 &= ( \vec{e}_1 , \cev{f}_1, \cev{e}_{3}, \vec{f}_1 \ ), \\ 
C_{2} &= ( \vec{e}_{2} , \cev{f}_{3}, \cev{e}_{1}, \vec{f}_{2}),\\  
C_{3} &= (\vec{e}_{3} , \cev{f}_{6}, \cev{e}_{2}, \vec{f}_{5}),  \\
C_{4} &= ( \cev{f}_{2} , \vec{f}_{4}, \vec{f}_{3}, \cev{f}_{4} ),  \\
C_{5} &= (  \cev{f}_{5} , \vec{f}_{7}, \vec{f}_{6}, \cev{f}_{7}).
\end{align*}
Now, we count the boundary components of $\Gamma_5$ using the permutation $\sigma_{\infty}$. The boundary components of the graph $\Gamma_5$ given by: \\
$\partial_1= \vec{e_1} \vec{f_2} \vec{f_4} \cev{f_2} \vec{e_2} \vec{f_5} \vec{f_7} \cev{f_5} \vec{e_3} \vec{f_1} \cev{e_3} \cev{f_6} \cev{f_7} \vec{f_6} \cev{e_2} \cev{f_3} \cev{f_4} \vec{f_3} \cev{e_1} \cev{f_1}.$\\
So, there is only one boundary component.
\section{Upper Bound On Mapping Class Orbits of Maximum Size Filling }
Let $ I(S_g) $ represent the set of all isotopy classes of simple closed curves on $S_g.$ $ I(S_g)^{n} $ denotes the Cartesian product of $n \in \mathbb{N}$ copies of $ I(S_g) $. Define an equivalence relation $\sim$  on $I(S_g)^n $ by $ (\gamma_1, \gamma_2, \ldots, \gamma_n) \sim (\gamma_{\sigma(1)}, \gamma_{\sigma(2)}, \ldots, \gamma_{\sigma(n)}) , \text{ where  $\sigma$  is a permutation of } \{1, 2, \ldots, n\}$. The mapping class group $Mod(S_g)$ \text{ acts on } $I(S_g)^n$ by $[f] \cdot (\gamma_1, \gamma_2, \ldots, \gamma_n) = (f \circ \gamma_1, \ldots, f \circ \gamma_n)$. This action induces an action on the quotient space $ \mathcal{I} (S_g):= I(S_g)^n / \sim$. Let $X=(c_1, \dots c_n) \in \mathcal{I} (S_g)$ be a filling of $S_g$ and $[f] \in Mod (S_g)$, then $[f].X=(f \circ c_1, \ldots, f \circ c_n)$ is also a filling of $S_g$ with same number of boundary components. In this section, we prove Theorem~\ref{orbit}.
\begin{proof}[ Proof of
Theorem \ref{orbit}] First we will consider when $b=1$ and $g\geq 2$. For $b=1$, the fat graph structure consists of 1 cycle of length $2 \leq s \leq g$, $s$ cycles of length 1 and $2g-1-s$ cycles of length 2. At the cycle of length $s$, one vertex has chain of odd length and the remaining $s-1$ vertices have chains of even length. Hence, we get the equation $2N_{1}+\sum_{i=2} ^{s} 2 N_i+1=2g-1-s$ which implies $\sum_{i=1} ^{s} N_i=g-s$. We define $\zeta^{1}_{s} =\{ (N_1, \cdots N_s) | N_{i} \in \mathbb{N} \cup \{0\}, \sum N_{i}=g-s\}/ \sim $ where the equivalence relation $\sim$ is defined as $ (N_{1}, \cdots , N_s)\sim(N^{'}_{1}, \cdots , N^{'}_s)$ if there exist $\sigma \ \in S_{s}$ such that $N^{'}_{i}=N_{\sigma_(i)} \text{ for } 1 \leq i \leq s$. Then define $\Re^{s}_{g-s}=\zeta^{1}_{2} \cup \zeta^{1}_{3}\cup \cdots \cup \zeta^{1}_{s} $. Hence, $N_{g,1}\geq |\Re^{s} _{g-s}|+1$.

 Next, we consider the case when $b=2$ and $g\geq 2$. For $b=2$, the fat graph structure consists of 1 cycle of length $2 \leq s \leq g+1$, $s$ cycles of length 1 and $2g-s$ cycles of length 2. At the cycle of length $s$,two vertices have chains of odd length and the the remaining $s-2$ vertices have chains of even length. Hence we get the equation $2N_{1}+2N_{2}+\sum_{i=3} ^{s} 2 N_i+1=2g-s$ which implies $\sum_{i=1} ^{s} N_i=g+1-s$. We define $\zeta^{2}_{s} =\{ (N_1, \cdots N_s) | N_{i} \in \mathbb{N} \cup \{0\}, \sum N_{i}=g+1-s\}/ \sim $ where equivalence relation $\sim$ defined as $ (N_{1}, \cdots , N_s)\sim(N^{'}_{1}, \cdots , N^{'}_s)$ if there exist $\sigma \ \in S_{s}$ such that $N^{'}_{i}=N_{\sigma_(i)} \text{ for } 1 \leq i \leq s$. Then define $\Re^{s}_{g-s+1}=\zeta^{2}_{2} \cup \zeta^{2}_{3}\cup \cdots \cup \zeta^{2}_{s} $. Hence, $N_{g,2}\geq |\Re^{s} _{g-s+1}|+1$.

 Finally, we consider the case when $3 \leq b \leq g-2$. The fat graph structure for this case consists of the $1$ cycle of length $2b$, $2g-b-2$ cycles of length $2$ and $2b$ cycles of length $1$. Hence, we get the equation: $\sum_{n=1} ^{b} 2 N_i-b=2g-2-b $ which implies $ \sum_{n=1} ^{b}  N_i=g-1 $. We define $\zeta^{b}_{b} =\{ (N_1, \cdots N_b) | N_{i} \in \mathbb{N} \cup \{0\}, \sum N_{i}=g-1\}/ \sim $ where equivalence relation $\sim$ defined as $ (N_{1}, \cdots , N_b)\sim(N^{'}_{1}, \cdots , N^{'}_b)$ if there exist $\sigma \ \in S_{b}$ such that $N^{'}_{i}=N_{\sigma_(i)} \text{ for } 1 \leq i \leq b$.Then define $\Re^{b}_{g-1}=\zeta^{b}_{3} \cup \zeta^{b}_{4}\cup \cdots \cup \zeta^{b}_{b} $. Hence, $N_{g,b}\geq |\Re^{b}_{g-1}|+1$.\\



\end{proof}

\section*{Acknowledgements}
The First author would like to thank Council of Scientific and Industrial Research (CSIR)( 09/1290(0002)2020-EMR-I ) for providing financial support. The second author gratefully acknowledges the financial support from the Science and Engineering Research Board (SERB), Government of India through MATRICS grant (File Number: MTR/2021/000067).

\bibliographystyle{alpha}
\bibliography{Reference}
\end{document}